\documentclass[11pt,reqno]{amsart}
\usepackage{euscript,verbatim,epsfig}

\newtheorem{theorem}{Theorem}
\newtheorem{proposition}[theorem]{Proposition}
\newtheorem{corollary}[theorem]{Corollary}

\newtheorem*{thma}{Theorem 5a}
\def\vol{{\rm vol} }
\newcommand{\beha}{\begin{enumerate}}
\newcommand{\behe}{\end{enumerate}}
\renewcommand{\epsilon}{\varepsilon}

\newcommand{\bR}{{\mathbb R}}
\newcommand{\bC}{{\mathbb C}}
\newcommand{\bZ}{{\mathbb Z}}
\newcommand{\bN}{{\mathbb N}}

\let\L=\Lambda
\let\e=\epsilon
\let\d=\delta

\DeclareMathSymbol{\varnothing}{\mathord}{AMSb}{"3F}

\begin{document}

\title{Stable sets, hyperbolicity and dimension}

\begin{abstract}
In this note we derive an upper bound for the Hausdorff  dimension of the
stable set of a hyperbolic set $\Lambda$ of a $C^2$
diffeomorphisms on a $n$-dimensional manifold. As a consequence we
obtain that $\dim_H W^s(\Lambda)=n$ is equivalent
to the existence of a SRB-measure. We also discuss related
results in the case of expanding maps.
\end{abstract}

\date{\today}

\author{Rasul Shafikov} \address{Department of Mathematics, SUNY at
  Stony Brook, Stony Brook, NY 11794}
\email{shafikov@math.sunysb.edu}

\author{Christian Wolf} \address{Department of Mathematics, Wichita
  State University, Wichita, KS 67260, USA} 
\email{cwolf@math.wichita.edu}

%\subjclass{Primary: 11K55, 37C45.}
\maketitle

%-------------------------------------
\section{Introduction}
Let $M$ be a $n$-dimensional smooth Riemannian manifold, $f:M\to M$ a
diffeomorphism, and $\Lambda\subset M$ a locally
maximal hyperbolic set of $f$. We define the stable set of $\Lambda$
by 
\[
W^s(\Lambda)=\{x\in M: {\rm dist}(f^k(x), \Lambda)\to 0\ \text{for}\
k\to\infty\}. 
\]     
In this paper we investigate the complexity of $W^s(\Lambda)$ in terms
of its Hausdorff dimension $\dim_H W^s(\Lambda)$. 
It is a classical result of Bowen \cite{Bo1} that there exist examples
of $C^1$ horseshoes $\Lambda$ with positive Lebesgue measure, in
particular $\dim_H W^s(\Lambda)=n$. On the other hand, by a result of   
Bowen and Ruelle \cite{BR},  if $f$ is a $C^2$-diffeomorphism and
$\Lambda$ is not an attractor then $W^s(\Lambda)$ has zero Lebesgue
measure. We extend the latter result by showing  that the Hausdorff
dimension of $W^s(\Lambda)$ is strictly smaller than $n$. More
precisely, we derive an upper bound for the Hausdorff dimension of
$W^s(\Lambda)$ which is given in terms of the exponential expansion rate
of the tangent map defined by
\begin{equation}\label{s}
s=\lim_{k\to\infty}\frac{1}{k}\log\left(\max\{||Df^k(x)||: x\in
\Lambda\}\right),   
\end{equation}
and the topological pressure of the unstable Jacobian. If
$\Lambda$ is not an attractor this bound is strictly smaller than $n$. 

\begin{theorem}\label{mainresult}
Let $f:M\to M$ be a $C^2$-diffeomorphism, and $\Lambda$  a  locally
maximal hyperbolic set of $f$, which is not a periodic orbit, such that 
$f|\Lambda$ is topologically mixing. Define $\phi^u=-\log |\det
Df|{E^u}|$. Then 
\begin{equation}\label{maininequalityint}
\dim_H W^s(\Lambda)\leq n+\frac{P(\phi^u)}{s}.
\end{equation}
\end{theorem}

Here $P(\phi^u)$  denotes  the topological pressure of $\phi^u$, see
Section~2 for details. We note that our result holds for manifolds of
arbitrary dimension, in particular, we do not require that $f$ is
conformal on $\Lambda$. Since $f$ is a diffeomorphism, the analogous 
result also holds for the unstable set $W^u(\Lambda)$ of $\Lambda$.

It follows from the proof of Theorem \ref{mainresult} that the upper
bound in \eqref{maininequalityint}  also provides an upper bound for
the upper box dimension of  certain subsets of $W^s(\Lambda)$. In
particular, the following holds. 

\begin{corollary}\label{corref1}
Suppose that $f$ and $\Lambda$ are as in Theorem \ref{mainresult}, and
assume that $\Lambda$ has  empty interior. Then the upper box dimension
$\overline{\dim}_B \Lambda$ of $\Lambda$ is strictly smaller than
$n$. 
\end{corollary}

We note that the case $\Lambda$  having  non-empty interior actually
occurs. For example, if 
$f$ is an Anosov diffeomorphism, then the entire manifold $M$ is a
hyperbolic set.  
The result of Corollary \ref{corref1} is known in some special
cases. In particular, if $M$ is a surface, that is, $n=2$, then the
classical result of McCluskey and Manning \cite{MM} states that the
Hausdorff dimension of $\Lambda$ coincides with  the box dimension of
$\Lambda$, and that its  value is strictly smaller than $2$. Recently,
significant progress has also been made  towards the estimation of the
Hausdorff and box dimension of  hyperbolic sets of higher dimensional
manifolds, see for instance \cite{Ba}, \cite{F}, and  the references
in \cite{P}. These estimates are typically given in terms of the
topological pressure of functions related to the expansion/contraction
rates of the tangent map on the hyperbolic splitting. However, it is
often difficult to calculate these upper bounds in particular
examples. Corollary \ref{corref1}, on the other hand, states without
further calculations that the dimension of a hyperbolic set is
strictly smaller than $n$. This has not been known before in this
generality.  

Another consequence of Theorem \ref{mainresult} is a  new
characterization  for $f$ to have an invariant probability measure
$\mu$ supported on $\Lambda$ whose conditional measures on the
unstable manifolds are absolutely continuous with respect to the
Lebesgue measure. Such a measure $\mu$ is called a SRB measure of the
diffeomorphism $f$. 
 
\begin{theorem}\label{mainresult2}
Suppose that $f$ and $\Lambda$ are as in Theorem \ref{mainresult}.
Then the following are
equivalent. 
\begin{enumerate}
\item[(i)]
$\dim_H W^s(\Lambda)=n$;
\item[(ii)]
$f$ admits a SRB measure on $\Lambda$;
\item[(iii)]
$\Lambda$ is an attractor of $f$. 
\item[(iv)]
$W^s(\Lambda)$ has positive Lebesgue measure
\end{enumerate}
\end{theorem}

The novelty of Theorem \ref{mainresult2} is the implication (i)
implies (ii), while the other implications are well-known, and are not
consequences  of our results.  

This paper is organized as follows. In Section 2 we consider
hyperbolic diffeomorphisms and derive an upper bound of the box
dimension of particular sets associated with a hyperbolic
set. Furthermore, we apply this bound to establish the results stated
in the introduction. Finally, in Section 3 we study repellers of
expanding maps, and derive related results as in the case of
diffeomorphisms.  

\section{Stable sets for diffeomorphisms}

Let $f\colon M\to M$ be a $C^{2}$ diffeomorphism on a
$n$-dimensional Riemannian manifold, and $\Lambda\subset M$ a 
hyperbolic set.  This means that $\Lambda$ is a compact $f$-invariant set, and that
there exist a continuous splitting of the tangent bundle 
$T_\Lambda M=E^u\oplus E^s$, and constants $c>0$ and $\lambda\in(0,1)$
such that for each $x\in\Lambda$:
\begin{enumerate}
\item $Df(x)(E^u_x)=E^u_{f(x)}$ and $Df(x)(E^s_x)=E^s_{f(x)}$;
\item $\lVert Df^{-k}(x) v\rVert\leq c\lambda^k\|v\|$ whenever $v\in
E^u_x$ and $k\in\bN$;
\item $\lVert Df^k(x)v\rVert\leq c\lambda^k\|v\|$ whenever $v\in
E^s_x$ and $k\in\bN$.
\end{enumerate}
We say that $\Lambda$ is locally maximal if
there exists an open neighborhood $U$ of $\Lambda$ such that
$\Lambda=\bigcap_{k\in\bZ}f^k U$.  We shall always assume that
$\Lambda$ is a locally maximal  hyperbolic set, which is not a
periodic orbit, and that $f|\Lambda$ is topologically mixing.  

We say that $\Lambda$ is an attractor of $f$ if there are arbitrarily
small neighborhoods $U$ of $\Lambda$ such that $f(U)\subset U$.
This includes the case of Anosov diffeomorphisms for
which the entire manifold $M$ is a hyperbolic set. We say that
$\Lambda$ is a repeller of $f$ if it is an attractor of $f^{-1}$. We
note that there are other definitions of an attractor and
repeller. Here we follow the notation of \cite{Bo2} which suits well
our purposes. 

Following Bowen \cite{Bo2}, for $x\in \Lambda$, $\epsilon>0$ and $k\in 
\bN$ we define the sets
\begin{equation}\label{eqbo1}
B(x,\epsilon,k)=\{y\in M: |f^i(x)-f^i(y)|<\epsilon, i=0,\dots,k-1\},
\end{equation}
where $|\cdot|$ denotes the distance induced by the Riemannian
metric. Define
\begin{equation}\label{last}
B(\Lambda,\epsilon,k)=\bigcup_{x\in \Lambda} B(x,\epsilon,k). 
\end{equation}
Furthermore, we define the local stable set of
$\Lambda$ by
\begin{equation}
W^s_\epsilon(\Lambda)=\{x\in M: {\rm
  dist}(f^k(x),\Lambda)<\epsilon\ {\rm for\ all}\ k\in \bN \}.
  \end{equation}
Analogously, we define the local unstable set $W^u_\epsilon(\Lambda)$
of $\Lambda$ by considering $f^{-1}$. It is an immediate  consequence
of the shadowing lemma that if $\epsilon$ is sufficiently small, then
the exists $\overline{\epsilon}>0$ such that 
\begin{equation}\label{e-bar}
W^s_\epsilon(\Lambda)\subset \bigcup_{x\in\Lambda}
W^s_{\overline{\epsilon}}(x),
\end{equation}
where $W^s_{\overline{\epsilon}}(x)$ denotes the local stable manifold
of size $\overline \epsilon$ of $x\in \Lambda$. Furthermore, by choice
of $\epsilon$, the number $\overline\epsilon$ can be chosen
arbitrarily small. We will also need the topological pressure 
$P(\varphi)$ of a continuous function $\varphi:\Lambda\to\bR$ (see
\cite{Wa} for the definition and details). We are interested in the
particular function  $\phi^u:\Lambda\to \bR$ defined  by 
 $\phi^{u}(x)=-\log \vert{\rm det}Df(x)|{E^u_x}\vert$, where ${\rm
  det}Df(x)|{E^u_x}$ denotes the Jacobian of the linear map
$Df(x)|E^{u}_x$. We will need the following result which is  due to
Bowen.

\begin{proposition}[\cite{Bo2}]\label{thbowen}
  If $\epsilon>0$ is small enough then:
  \begin{enumerate}
  \item[(i)] $
    \lim_{k\to\infty}\frac{1}{k}\log\left(\vol(B(\Lambda,\epsilon,k)\right)=
    P(\phi^u)\leq 0;$ 
  \item[(ii)]  $P(\phi^u)=0$ if and only if $W^s_\epsilon(\Lambda)$
  has nonempty interior, in which case $\Lambda$ is an attractor. 
  \end{enumerate}
\end{proposition}

We note that the right-hand side inequality in (i) also follows as an
application of the Margulis-Ruelle inequality and the variational
principle. The following is the main result of this paper.

\begin{theorem}\label{thmainbox}
  Let $f:M\to M$  be a $C^2$ diffeomorphism on a $n$-dimensional
  Riemannian manifold, and let $\Lambda\subset M$ be a locally
  maximal hyperbolic set, which is not a periodic orbit, such that
  $f|\Lambda$ is topologically mixing. Then for sufficiently small
  $\epsilon>0$ we have  
  \begin{equation}\label{maininequality}
    \overline{\dim}_B W^s_\epsilon(\Lambda)\leq
    n+\frac{P(\phi^u)}{s},
  \end{equation}
  where $s$ is defined as in \eqref{s}.
\end{theorem}

\noindent{\it Remark.} It follows from Proposition \ref{thbowen} that
if $\Lambda$ is not an attractor, then $P(\phi^u)<0$, and therefore
inequality \eqref{maininequality} provides a non-trivial estimate.

\begin{proof}[Proof of Theorem \ref{thmainbox}.]
First observe that since the operator norm is submultiplicative, the
limit defining $s$ exists (see e.g. \cite{Wa}). Furthermore, since
$\Lambda$ is not an attracting cycle, $s>0$.  If $\Lambda$ is an
attractor, then by Proposition~\ref{thbowen}, $P(\phi^u)=0$, and
inequality \eqref{maininequality} trivially holds. Thus, we may assume
that $\Lambda$ is not an attractor, in which case $P(\phi^u)<0$.

Let $\delta>0$. It follows from a simple continuity argument that
there exist $\epsilon>0$ and $k_\delta\in \bN$ such that for all 
$x\in B(W^{s}_\epsilon (\Lambda),\epsilon) =\{x\in M:\exists 
y\in W^{s}_\epsilon (\Lambda),\  |x-y|<\epsilon\} $ we have
\begin{equation}
\lVert Df^{k_\delta}(x)\rVert < \exp(k_\delta(s+\delta)).
\end{equation}
From now on we consider the map $g=f^{k_\delta}$. Note that $\Lambda$
is also a hyperbolic set of $g$. Evidently $W^{s}_\epsilon (\Lambda)$  
is forward invariant under $g$. It follows from the variational
principle that $P_g(\phi^u_g)=k_\delta P_f(\phi^u_f)$; moreover 
$s_g=k_\delta s_f$. Thus it is sufficient to prove inequality
\eqref{maininequality} for $g$. We continue to use the notation
$s,\phi^u, P(\phi^u),$ etc. for $g$ instead of $f$. 
Let $x\in \Lambda$ and
$k\in\bN$. It follows from Proposition \ref{thbowen} that for
$\epsilon$ sufficiently small, 
\begin{equation}\label{eqbo2}
P(\phi^u)=\lim_{k\to\infty}\frac{1}{k}\log(
\vol(B(\Lambda,2\epsilon,k))),
\end{equation} 
where $B(\Lambda,2\epsilon,,k)=\bigcup_{x\in\Lambda} B(x,2\epsilon,k)$ (see \eqref{eqbo1}). From
this we obtain that if $k$ is sufficiently large, then 
\begin{equation}\label{eqbo4}
\vol(B(\Lambda,2\epsilon,k))< \exp(k(P(\phi^u)+\delta)).
\end{equation}
For all $k\in \bN$ we define real numbers
\begin{displaymath}
r_k = \frac{\epsilon}{\exp(s+\delta)^k}
\end{displaymath}
and neighborhoods $B_k = B(W^{s}_\epsilon (\Lambda),r_k)$ of
$W^{s}_\epsilon (\Lambda)$. Let $y\in B_k$. Then there exists $ x\in
W^{s}_\epsilon (\Lambda)$ with $|x-y| < r_k$. An elementary induction
argument in combination with the mean-value theorem implies
$|g^{i}(x)- g^{i}(y)| < \epsilon$ for all $i\in\{0,\dots,k-1\}$.
Using \eqref{e-bar} and making $\epsilon$ smaller if necessary, we
can assure that $x$ is contained in the local stable manifold of size
$\epsilon$ of a point in $\Lambda$. It follows that $y\in
B(\Lambda,2\epsilon,k)$. Hence $B_k\subset B(\Lambda,2\epsilon,k)$.
Therefore, \eqref{eqbo4} implies that
\begin{equation}
\vol(B_k)<\exp(k(P(\phi^u)+\delta))
\end{equation}
for sufficiently large $k$.
Let us recall that for $t\in [0,n]$ the  $t$-dimensional upper
Minkowski content of a relatively compact set $A\subset M$ is defined
by 
\begin{displaymath}
M^{*t}(A)= \limsup_{\rho\to 0} \frac{\vol(A_\rho)}{(2\rho)^{n-t}},
\end{displaymath}
where $A_\rho=\{p\in M:\exists q\in A: |p-q|\leq \rho\}$, and $\vol$ denotes the volume induced by the Riemannian metric on $M$. Let
$t\in [0,n]$ and $\rho_k={\textstyle\frac{r_k}{2}}$ for all
$k\in\bN.$ Then we have
\begin{equation}\label{glej-4}
\begin{array}{lcl}
M^{*t} (W^s_\epsilon(\Lambda))&=&  \limsup\limits_{\rho\to 0}
\frac{\textstyle
\vol(W^s_\epsilon(\Lambda)_\rho)}{\textstyle(2\rho)^{n-t}}\\ &\leq&
\limsup\limits_{k\to\infty } \frac{\textstyle
\vol(W^s_\epsilon(\Lambda)_{\rho_k})}{\textstyle(2\rho_{k+1})^{n-t}}\\
&\leq &
 \limsup\limits_{k\to \infty} \frac{\textstyle \vol(B_k)}{\textstyle (r_{k+1})^{n-t}}\\
&\leq & \frac{\textstyle
\exp(s+\delta)^{n-t}}{\textstyle
\epsilon^{n-t}}\lim\limits_{k\to\infty}\left(\textstyle\exp(
s+\delta)^{n-t}\textstyle \exp(P(\phi^u)+\delta)\right)^k.
\end{array}
\end{equation}
Let  $t > n + \frac{P(\phi^u)+\delta}{ s+\delta}$.  Then
$\exp(s+\delta)^{n-t} \exp(P(\phi^u)+\delta) < 1.$ This implies
$M^{*t}(W^s_\epsilon(\Lambda)) = 0$, in particular, $t\geq
\overline{\dim}_B W^s_\epsilon(\Lambda)$. Since $\delta$
can be chosen arbitrarily small, the result follows.
\end{proof}

\noindent{\it Remark.} We note that the idea of estimating the
dimension of an invariant set by calculating  the volume of
neighborhoods of the set was introduced by the second author of this
paper in \cite{Wo} for estimating the dimension of an invariant set of
a $C^1$ diffeomorphism. For a related result in the case of polynomial 
automorphisms of $\bC^n$ see \cite{SW}. 

\bigskip

Applying Theorem \ref{mainresult} to $f^{-1}$ we obtain the analogous 
result for the local unstable set 
\[W^u_\epsilon(\Lambda)=\{x\in M: {\rm dist}(f^{-k}(x), \Lambda)<
\epsilon \ \ \text{for all}\ k\in \bN\}.\]
     
\begin{thma}\label{thupchris}
  Let $f:M\to M$  be a $C^2$ diffeomorphism on a $n$-dimensional
  Riemannian manifold, $\Lambda\subset M$ a locally maximal
  hyperbolic set, which is not a periodic orbit, such that $f|\Lambda$
  is topologically mixing. Let 
  $\phi^s=\log |\det Df|E^s|$. Then for sufficiently small
  $\epsilon>0$  we have  
  \begin{equation}\label{maininequality22}
    \overline{\dim}_B W^u_\epsilon(\Lambda)\leq n+\frac{P(\phi^s)}{s},
  \end{equation}
  where $s$ is defined as in \eqref{s} for the map $f^{-1}$.
\end{thma}

Note that if $M$ is a surface, then the dimension of
$W^{u/s}(\Lambda)$ can be expressed using Bowen's formula. Namely, 
\[
t^{u/s}:=\dim_H W^{u/s}(x)\cap \Lambda= \overline{\dim}_B
 W^{u/s}(x)\cap \Lambda
\] 
are independent of $x\in \Lambda$ (see \cite{MM}), and it is not to
hard to see that
\begin{equation}\label{13}
\dim_H W^{s/u}_\epsilon(\Lambda)=\overline{\dim}_B
W^{s/u}_\epsilon(\Lambda)= t^{u/s}+1 . 
\end{equation}
The right equality in \eqref{13} can be shown using the fact that when
$n=2$ the holonomies are lipschitz continuous (see for instance
\cite{KH}). On the other hand, the following example shows that the
dimension of $W^s_\epsilon(\Lambda)$ can be arbitrarily close to $n$. 

\bigskip

\noindent
{\bf Example 1. }
Let $B\subset \bR^2$ be a unit square, and let $f:B\to\bR^2$ be a
linear horseshoe  map with the expansion rate $\lambda^u>2$ and the
contraction rate $\lambda^s<\frac{1}{2}$, see Figure 1. It is
well-known that $\Lambda=\{x\in B: f^k(x)\in B\ \text{for all}\ k\in
\bZ\}$ is a locally maximal hyperbolic set of $f$, and that
$f|\Lambda$ is topologically mixing. Moreover,  we have  
\[
t^u=\frac{\log 2}{\log \lambda^u}\ \text{\ \ and\ \ }\ t^s=-\frac{\log 2}{\log
  \lambda^s}. 
\]

\begin{picture}(400,220)(10,20)
%% this is outer frame
%\put(0,0){\framebox(400,240)}
%% this is the title
\put(125,4){{\bf Figure 1.} A Linear Horseshoe Map}
%% this is the first square
\put(23,90){\framebox(80,80)}
%% first square ones
\put(10,127){$1$}
\put(60,76){$1$}
%% this is the third square
\put(270,90){\framebox(80,80)}
%% this is a rectangle
\put(175,40){\framebox(16,190)}
%% the first arrow
\put(120,130){\vector(1,0){25}}
%% this is the second arrow
\put(223,130){\vector(1,0){25}}
%% four vertical lines of the horseshoe
\put(286,70){\line(0,1){100}}
\put(302,70){\line(0,1){100}}
\put(318,70){\line(0,1){100}}
\put(334,70){\line(0,1){100}}
%% bottom of the horseshoe
\put(286,70){\line(1,0){16}}
\put(318,70){\line(1,0){16}}
%% the upper hat of the horseshoe
\put(310,170){\oval(48,48)[t]}
%% the lower hat of the horseshoe
\put(310,170){\oval(16,16)[t]}
%% first lambda
\put(155,125){$\lambda^u$}
%% second lambda
\put(179,26){$\lambda^s$}
%% third and fourth lambda
\put(290,55){$\lambda^s$}
\put(322,55){$\lambda^s$}
\end{picture}

\bigskip

\bigskip

As it was mentioned above, we also have $\dim_H
W^{s/u}_\epsilon(\Lambda)=\overline{\dim}_B
W^{s/u}_\epsilon(\Lambda)=t^{u/s}+1$. Therefore, by choosing
$\lambda^u$ close to $2$ (respectively $\lambda^s$ close to
$\frac{1}{2}$) we obtain the following. 

\begin{corollary}
For each $\epsilon>0$ there exists a linear horseshoe map of $\bR^2$ 
such that $\dim_H W^{s/u}_\epsilon(\Lambda)=\overline{\dim}_B W^{s/u}_\epsilon(\Lambda) > 2-\epsilon$.  
\end{corollary}

We now provide the proofs of the results stated in the introduction. 

\begin{proof}[Proof of Theorem \ref{mainresult}] 
  It follows from Theorem  
  \ref{thmainbox} that $\dim_H W^s_\epsilon(\Lambda)\leq
  n+\frac{P(\phi^u)}{s}$. Obviously, $W^s_\epsilon(x)\subset
  W^s_\epsilon(\Lambda)$ for all $x\in\Lambda$, hence 
  \begin{equation}\label{eqlala}
    \dim_H \bigcup_{x\in\Lambda} W^s_\epsilon(x)\leq
    n+\frac{P(\phi^u)}{s}. 
  \end{equation}
  It is a consequence of the shadowing lemma that
  \begin{equation}
    W^s(\Lambda)=\bigcup_{x\in\Lambda} W^s(x).
  \end{equation}
  On the other hand, we have $W^s(x)=\bigcup_{k\in \bN}
  f^{-k}(W^s_\epsilon(f^k(x)))$. Together we obtain 
  \begin{equation}\label{eqhaha}
    W^s(\Lambda)=\bigcup_{k\in\bN} f^{-k}\left(\bigcup_{x\in\Lambda}
    W^s_\epsilon(x)\right). 
  \end{equation}
  The theorem follows now from equations \eqref{eqlala}, \eqref{eqhaha} and the
  fact that the Hausdorff dimension is stable with respect to countable unions. 
\end{proof}

\noindent{\it Remark.} Since we have taken countable unions of sets
whose box dimensions are uniformly bounded above by
$n+\frac{P(\phi^u)}{s}$, we obtain the same upper bound for the
  packing dimension of $W^s(\Lambda)$.

\begin{proof}[Proof of Corollary \ref{corref1}] If $\Lambda$ is an
  attractor, then $\Lambda$ can not be a 
  repeller. Otherwise, we would have $W^{u/s}_\epsilon(x)\subset
  \Lambda$ for all $x\in\Lambda$, in which case, since $\Lambda$ has a
  local product structure, $\Lambda$ would have a non-empty
  interior. Therefore, Proposition
  \ref{thbowen} (or the analogous proposition for $f^{-1}$) implies
  that either $P(\phi^u)<0$ or $P(\phi^s)<0$. The result now follows
  from  Theorems \ref{thmainbox} and 5a. 
\end{proof}

\begin{proof}[Proof of Theorem \ref{mainresult2}]
  (i)$\Rightarrow$(ii) Assume $\dim_H W^s(\Lambda)=n$. Then it follows
  from Theorem \ref{mainresult}  and Proposition \ref{thbowen} (i) that
  $P(\phi^u)=0$. Since $\Lambda$ is a locally maximal hyperbolic set
  on which $f$ is topologically mixing, there exists a unique
  equilibrium measure $\mu_{\phi^u}$ of the potential $\phi^u$. This
  means that 
  \begin{equation}\label{eqrgh} 
    0=P(\phi^u)=h_{\mu_{\phi^u}}(f) + \int \phi^u d\mu_{\phi^u},
  \end{equation} 
  where $h_{\mu_{\phi^u}}(f)$ denotes the measure-theoretic entropy of
  $f$ with respect to $\mu_{\phi^u}$. Moreover, $\mu_{\phi^u}$ is
  ergodic. Let $\lambda_1(x)<\cdots<\lambda_l(x)$ be the Lyapunov exponents
  of $x$ with respect to $f$ with multiplicities
  $m_1(x),\cdots,m_l(x)$. The fact that $\mu_{\phi^u}$ is ergodic
  implies that the Lyapunov exponents and the
  multiplicities are constant
  $\mu_{\phi^u}$-almost everywhere.  We denote the corresponding
  values by $\lambda_i(\mu_{\phi^u})$ and $m_i(\mu_{\phi^u})$. It
  follows from the fact 
  that $\Lambda$ is a hyperbolic set  that 
  \begin{equation}\label{eq1678}
    -\int \phi^u d\mu_{\phi^u}=\sum_{\lambda_i(\mu_{\phi^u})>0}
     \lambda_i(\mu_{\phi^u}) m_i(\mu_{\phi^u}).
  \end{equation}
  Therefore, equations \eqref{eqrgh} and \eqref{eq1678} imply that  
  \begin{equation}
    h_{\mu_{\phi^u}}(f)=\sum_{\lambda_i(\mu_{\phi^u})>0}
    \lambda_i(\mu_{\phi^u})  m_i(\mu_{\phi^u})=
    \int\sum_{\lambda_i(x)>0} \lambda_i(x) m_i(x)\ d\mu_{\phi^u}. 
  \end{equation}
  Hence $\mu_{\phi^u}$ satisfies Pesin's entropy formula. Finally,
  from the work of Ledrappier, Strelcyn and Young in \cite{LS} and
  \cite{LY} we conclude that $\mu_{\phi^u}$ is a SRB measure. 
  
  (ii)$\Rightarrow$(iii) follows from \cite{Bo2} and
  \cite{LY}. Finally, (iii)$\Rightarrow$ (iv) and (iv)$\Rightarrow$(i)
  are trivial.
\end{proof}

\section{Expanding maps} 

In this section we consider a  $C^2$ map $f$
on a $n$-dimensional smooth Riemannian manifold $M$ to itself. The
map $f$ is not assumed to be invertible. Let $\Lambda\subset M$ be a
compact invariant set of $f$. We say that $f$ is expanding on
$\Lambda$  if there exist $c>0$ and $\lambda\in (1,\infty)$ such
that for each $x\in \Lambda$: 
\[
\lVert Df^k(x) v\rVert\geq c\lambda^k\|v\|\quad {\rm
  whenever}\quad v\in T_xM\ {\rm and}\ k\in\bN.
\]
Furthermore, we say that $\Lambda$ is locally maximal if there exists
an open neighborhood $U\subset M$ of $\Lambda$ such that
$\Lambda=\bigcap_{k\in\bN} f^{-k}(U)$. If $f$ is expanding on a
locally maximal set $\Lambda$  we say that $\Lambda$ is a repeller of
$f$. We shall always assume that $f$ is expanding on $\Lambda$,
$\Lambda$ is locally maximal, and $f|\Lambda$ is topologically
mixing. We define the function $\phi:\Lambda\to\bR$ by
$\phi=-\log|\det Df|$.  

 We start by proving the version of Proposition \ref{thbowen} for
 expanding maps. 

\begin{proposition}\label{propbowen1}
  If $\epsilon>0$ is small enough then
  \begin{equation}\label{eqrcd}
    \lim_{k\to\infty}\frac{1}{k}\log\left(\vol(B(\Lambda,\epsilon,k)\right)\leq
    P(\phi)\leq 0
  \end{equation}
\end{proposition}

\begin{proof} We first show the right-hand side inequality of
\eqref{eqrcd}. Consider an ergodic invariant probability measure $\mu$
supported on $\Lambda$, and let $0<\lambda_1(\mu)<\cdots
<\lambda_l(\mu)$ be the  Lyaponov exponents of $\mu$. Denote by
$m_i(\mu)$ the multiplicity of the Lyapunov exponent
$\lambda_i(\mu)$. It follows from the Margulis-Ruelle inequality that   
\begin{equation}\label{eqwer1}
h_\mu(f)\leq \sum_{i} \lambda_i(\mu)\ m_i(\mu) = 
-\int \phi d\mu.
\end{equation}
On the other hand, the variational principle states that 
\begin{equation}\label{eqwer2}
P(\phi)=\sup_\mu\left(h_\mu(f)+\int \phi d\mu\right),
\end{equation}
where the supremum is taken over all (ergodic) invariant probability
measures on $\Lambda$. Combining \eqref{eqwer1} and \eqref{eqwer2}
yields $P(\phi)\leq 0$. 

%mystuff
We now show the left-hand side inequality in \eqref{eqrcd}.
Fix some $\d\le \e$. We say that  a set $E\subset M$ is
$(k,\delta)$-separated if for all $y, z\in E$, $y\not=z$, there exists
$i\in \{0,...,k-1\}$ with $|f^i(y)- f^i(z)|\ge \delta$. Let
$E_k(\delta)$ be a maximal $(k,\delta)$-separated subset of $\Lambda$. 
Let  $x\in \Lambda$; then
$x\in B(y,\d,k)$ for some $y\in E_k(\delta)$, because 
otherwise $E_k(\d)\cup\{x\}$ would be $(k,\delta)$-separated. Here
$B(y,\d,k)$ is defined as in \eqref{eqbo1}. We conclude that
$B(x,\e,k)\subset B(y,\d+\e,k)$, and 
\begin{equation}\label{verylast}
B(\L,\e,k)\subset \bigcup_{y\in E_k(\d)}B(y,\d+\e,k),
\end{equation}
see \eqref{last} for the definition.
For $x\in \Lambda$ and $k\in\bN$ we define
\[
S_k\phi(x)=\sum_{i=0}^{k-1}\phi(f^i(x)).
\]

Analogously as in  the case of diffeomorphisms (see
\cite{Bo2}, \cite{BR}), there exists $C_{\delta+\epsilon} >1$ such
that if $\epsilon$ is small enough then 
\begin{equation}\label{veryverylast}
\vol(B(x,\delta+\epsilon,k))\leq C_{\delta+\epsilon}\cdot
\exp S_k \phi(x)
\end{equation}
for all $x\in\Lambda$ and $k\in \bN$. Therefore \eqref{verylast} and
\eqref{veryverylast} imply that
\begin{equation}\label{b1}
  \vol (B(\L,\e,k))\le C_{\d+\e}\cdot\sum_{y\in E_k(\d)} \exp
  S_k\phi(y).
\end{equation}
 The map $f|\Lambda$ is
expansive; therefore for $\delta$ small enough we have
\begin{equation}\label{eqwert}
P(\phi)=\limsup_{k\to\infty}\frac{1}{k}\log \left( \sum_{y\in
    E_k(\delta)}\exp (S_k\phi(y))\right).
\end{equation}
The result follows now by taking the logarithm, dividing by $k$ and
taking the upper limit in \eqref{b1} and applying \eqref{eqwert}. 
%%end of mystuff
\end{proof}

We would like to point out that the proof of Proposition
\ref{propbowen1} is based on the ideas of the corresponding proof of
Bowen \cite{Bo2} for diffeomorphisms. 

Quan and Zhu classified in \cite{QZ} the invariant  measures $\mu$
of a $C^2$ endomorphism $f$ for which Pesin's entropy formula holds;
they showed that $\mu$ has this property if and only if $\mu$
satisfies the SRB property. This property is a generalized condition
of an SRB measure for diffeomorphisms defined for the corresponding
measure on the inverse limit map, see \cite{QZ} for details. In
particular, if $f$ is an expanding map and $\mu$ is absolutely
continuous with respect to Lebesgue, then Pesin's entropy formula
holds \cite{H}. 

\begin{proposition}\label{propbowen2}
$P(\phi)= 0$ if and only if $f$ admits an invariant measure $\mu$
satisfying  the SRB-property. In particular this is the case when $f$
has an invariant measure which is absolutely continuous with respect
to Lebesgue. 
\end{proposition}

\begin{proof}
If $P(\phi)=0$ then the same arguments as in the proof of Theorem
\ref{mainresult2} imply that there exists an invariant measure $\mu$
satisfying Pesin's formula; hence, $\mu$ has the SRB property. On the
other hand, if $\mu$ has the SRB property then $\mu$ satisfies Pesin's
entropy formula. It now follows from the variational principle and a
similar argument as in the proof of Theorem \ref{mainresult2} that
$P(\phi)=0$. 
\end{proof}

\noindent{\it Remark. } It is easy to see that $P(\phi)=0$ actually occurs. 
For example, if $f:S^1\to S^1, f(z)=z^2$, then the 
entire manifold $S^1$ is a repeller. In this case the measure of
maximal entropy (given by the distribution of the periodic
points of $f$) is absolutely continuous with respect to Lebesgue.  

\bigskip

We now present our main result for expanding maps.

\begin{theorem}\label{mainexpanding}
Let $f$ be a  $C^2$ self map  on $n$-dimensional smooth manifold $M$,
and let $\Lambda$ be a locally maximal repeller of $f$ such that
$f|\Lambda$ is topologically mixing. Then  
\begin{equation}\label{maininequalityendomorphisms}
\overline{\dim}_B \Lambda\leq n+\frac{P(\phi)}{s},
\end{equation}
where $s$ is defined as in \eqref{s}.
\end{theorem}

\begin{proof}
  The proof of the theorem is  analogous to the proof of Theorem
  \ref{thmainbox} just by replacing $W^s_{\epsilon}(\Lambda)$ by
  $\Lambda$ and applying  Proposition \ref{propbowen1}
  instead of Proposition~\ref{thbowen}.\end{proof}

As a consequence of Theorem \ref{mainexpanding} and Proposition
\ref{propbowen2} we obtain the following.

\begin{corollary}\label{mainresult2end}
  Let $f$ be a  $C^2$ self map of an $n$-dimensional smooth manifold
  $M$, and let $\Lambda$ be a repeller of $f$ such that $f|\Lambda$ is
  topologically mixing. If $\overline{\dim}_B \Lambda=n$ then
  $f$ admits an invariant measure $\mu$ satisfying the SRB property. 
\end{corollary}

\end{document}